\documentclass [a4paper,12pt]{article}
\usepackage[english]{babel}
\usepackage{amssymb,latexsym}
\usepackage{amsmath,amscd}                                                                                                           
\usepackage{mathrsfs}
\usepackage[all]{xy}
\usepackage{color}
\usepackage{epsfig,afterpage}
\usepackage{graphicx}
\usepackage{caption}
\usepackage{subcaption}

\usepackage[T1]{fontenc}
\usepackage{titling}
\usepackage{amsmath}
\usepackage{amssymb}
\usepackage{amsfonts}
\usepackage{multirow}
\usepackage{pifont}

\newtheorem{prop}{Proposition}[section]
\newtheorem{theo}[prop]{Theorem}
\newtheorem{Mtheo}[prop]{Main Theorem}

\newtheorem{remark}[prop]{Remark}
\newtheorem{defn}[prop]{Definition}

\newcommand{\gm}{\mbox{$\gamma$}}
\newcommand{\R}{\mbox{$\mathbb{R}$}}
 \newcommand{\bx}{\mbox{\boldmath $x$}}

\newcommand{\ba}{\mbox{\boldmath $a$}}

\newcommand{\bv}{\mbox{\boldmath $v$}}

\newcommand{\by}{\mbox{\boldmath $y$}}

\newcommand{\bu}{\mbox{\boldmath $u$}}

\newcommand{\bt}{\mbox{\boldmath $t$}}

\newcommand{\bn}{\mbox{\boldmath $n$}}

\date{\vspace{-5ex}}
\author{Graham Reeve}

\begin{document}

\title{Minkowski symmetry sets for 1-parameter families of plane curves}

\maketitle

\begin{abstract}
In this paper the generic bifurcations of the Minkowski symmetry set for 1-parameter families of plane curves are classified and the necessary and sufficient geometric criteria for each type are given.  The Minkowski symmetry set is an analogue of the standard Euclidean symmetry set, and is defined to be the locus of centres of all its bitangent pseudo-circles. It is shown that the list of possible bifurcation types are different to those that occur in the list
of possible types for the Euclidean symmetry set.      
\end{abstract}

\section{Introduction}
Symmetry sets and related constructions have provided useful representations of shapes for object recognition as well as attracting interest in their own right and in the geometric properties of curves that they reveal.
In the standard Euclidean plane, the {\it (Euclidean) symmetry set} of a curve $\gamma$ is defined as the locus of the centres
of circles that are tangent to $\gamma$ in at least two distinct points (bitangent), see for example \cite{BG86, bgg}.  The medial axis of $\gamma$ is a subset of its symmetry set, and is defined to be the locus of the centres
of circles that are bitangent to $\gamma$ and completely contained in $\gamma$. Introduced by Blum in 1967 \cite{blum}, the medial axis was originally designed as a tool for biological shape recognition and is also referred to as the central set, the topological skeleton, and the shock set for grassfire flows, see for example \cite{GibBras,KSMP}. 

The Minkowski symmetry set of a curve $\gamma$ was introduced in \cite{ST} as a Minkowski analogue of the (Euclidean) symmetry set. 
It is defined to be the locus of the centres of pseudo-circles that are bitangent to $\gamma$. In \cite{RT} the generic singularities of the Minkowski symmetry set were classified and the main result of the present paper is to extend this result to classify the generic singularities that occur for 1-parameter families of plane curves.  A well-studied, closely related construct to the symmetry set called the medial axis, is the locus of bitangent circles completely contained within $\gamma$ and has found various applications in computer vision. A Minkowski version of the medial axis was introduced in \cite{RT2}.

In \cite{BG86}, the transitions that occur for (Euclidean) symmetry sets of 1-parameter families of curves are classified. Moreover, the complete list of full bifurcation sets for a generic family of functions are given, and it is demonstrated that certain transitions are excluded for geometrical reasons. Analogous to this, in the present paper the generic bifurcations of the Minkowski symmetry set for 1-parameter families of plane curves are classified and their criteria are determined.



\begin{Mtheo}
The possible transitions types of the Minkowski Symmetry set for a generic curve are $A_1^4(a), A_1^4(b), A_2^2(a), A_2^2(b)$, $A_1A_3(a)$, $A_1A_3(b)$, $A_1^2A_2(a)$, $A_1^2A_2(b)$ and $A_4$.  
\end{Mtheo}


\begin{remark}
Note that the list of possible transitions types for the Minkowski Symmetry Set differs from that of the Euclidean Symmetry Set where only types $A_1^4(a), A_2^2(a), A_2^2(b), A_1A_3(a), A_1^2A_2(a)$ and $A_4$ can occur (see \cite{BG86}). 
\end{remark}

\begin{remark}
For the Euclidean Medial Axis, only types $A_1^4(a), A_1A_3(a)$ can the centre be on the medial axis (see for example \cite{GibKim}).  The Minkowski Medial Axis was defined in \cite{RT2} as the locus of centres of pseudo-circles that are bitangent to $\gamma$ with one of its branches.  It follows that only $A_1^4(b)$ and $A_1A_3(a)$ can have their centres on the Minkowski Medial Axis (see main text and Table on page \pageref{Table}). 
\end{remark}

\begin{remark}
In \cite{HoltTh} the an Affine version of the Symmetry Set called the Affine Distance Symmetry Set was considered.  It was shown that $A_1^4(a)$, $A_1^4(b)$,$ A_2^2(a)$, $A_2^2(b)$, $A_1A_3(a)$, $A_1A_3(b)$, $A_1^2A_2(a)$, $A_1^2A_2(b)$ and $A_4$ could occur generically.  In the case where $\gamma$ is an oval (a strictly convex, smooth and closed curve), it was also shown that $A_1^4(b),  A_1^2A_2(b)$ and $A_1A_3(b)$ were prohibited.
\end{remark}

\begin{center}
\label{Table}
\begin{tabular}{|c|c|c|c|}
  \hline
  &   {\rm \bf Euclidean} &  {\rm \bf Minkowski}  & {\rm \bf Affine}  \\
		\hline
 $A_1^4(a)$ & $\checkmark$ & Odd $\#$ points per branch & \checkmark \\
  \hline 
	$A_1^4(b)$ & $\times$ & Even $\#$ points per branch  &   Not for Ovals \\
  \hline 
   $A_2^2(a)$  & $\kappa_1'\kappa_2' > 0$  & $\kappa_1'\kappa_2' > 0$  (M. Curvature)  &  \checkmark \\
  \hline
	$A_2^2(b)$  & $\kappa_1'\kappa_2' < 0$  &   $\kappa_1'\kappa_2' < 0$  (M. Curvature) &  \checkmark \\
  \hline $A_1A_3(a)$ & \checkmark    & Points on different branches    & \checkmark  \\
  \hline $A_1A_3(b)$   &   $\times$  & Points on the same branch   &  Not for Ovals \\
  \hline $A_1^2A_2(a)$ & \checkmark  & $A_1$ points on same branch &   \checkmark \\
  \hline $A_1^2A_2(b)$  & $\times$   & $A_1$ points on opposite branches   &  Not for Ovals \\
  \hline $A_4$   &  \checkmark & \checkmark & \checkmark \\
  \hline
	\end{tabular}
	\end{center}

\section{The Minkowski pseudo-metric}\label{sec:prel}
The {\it Minkowski plane} $(\mathbb R_1^2, \langle, \rangle)$ is the vector space $\mathbb R^2$ endowed with the
pseudo-scalar product $\langle \bu, \bv \rangle = -u_0 v_0 + u_1v_1$, for any
$\bu = (u_0, u_1)$ and  $\bv = (v_0, v_1)$.
A vector $\bu \in \mathbb R^2_1$ is called {\it timelike} if $\langle \bu, \bu \rangle <0$, {\it spacelike} if $\langle \bu, \bu \rangle >0$, and
{\it lightlike} if $\langle \bu, \bu \rangle =0$.

 The norm of $\bu$ is defined by $||\bu|| = \sqrt{| \langle \bu, \bu \rangle |}$, and the perpendicular operator $\perp$ assigns ${\bf u}^\perp = (u_1, u_0)$.

There are three distinct types of pseudo-circles in $\mathbb R^2_1$ with centre $c \in \mathbb R^2_1$ and radius $r$, $r>0$, are defined as follows:
\begin{eqnarray}
\nonumber H^{1}(c, -r) &=& \{p \in \mathbb R^2_1\, | \, \langle p-c, p-c \rangle = -r^2 \}, \\
\nonumber S^{1}_1(c, r) &=& \{p \in \mathbb R^2_1 \, | \,  \langle p-c, p-c \rangle = r^2 \}, \\
\nonumber LC^*(c) &=& \{p \in \mathbb R^2_1 \setminus \{{\bf 0}\} \, | \,  \langle p-c, p-c \rangle = 0 \}.
\end{eqnarray}

Observe that $LC^*(c)$ is the union of the two lines through $c$ with tangent directions $(1,1)$ and $(1,-1)$, with the point
$c$ removed. The pseudo-circle $H^{1}(c, -r)$ has two branches which can be parametrised by $c+(\pm r\cosh(t),r\sinh(t))$, $t\in \mathbb R$.
The pseudo-circle $S^{1}(c, r)$ is also composed of two branches and these can be parametrised by $c+(r\sinh(t),\pm r\cosh(t))$, $t\in \mathbb R$.


Let $\gamma:S^1\to\mathbb R^2_1$ be an immersion, where $S^1$ is the unit Euclidean circle.
Call the curve $\gamma$ the image of the map $\gamma$ and say that it is a closed smooth curve (that is, $\gamma$
is a regular closed curve and may have points of self-intersection).

The curve $\gamma$ at $t_0$ is said to be spacelike if $\gamma '(t_0)$ is spacelike and is said to be timelike if $\gamma '(t_0)$ timelike.  
These are open properties so there is a neighbourhood of $t_0$ where the curve is either spacelike or timelike.
If $\gamma '(t_0)$ is lightlike then $\gamma(t_0)$ is said to be a lightlike point.
It is shown in \cite{ST} that the set of lightlike points of $\gamma$ is the union of at
least four disjoint non-empty and closed subsets of $\gamma$. 
The complement of these sets are disjoint connected spacelike or timelike pieces of the curve $\gamma$.

The spacelike and timelike components of $\gamma$ can be parametrised by arc length. Suppose that $\gamma(s)$, $s\in (\lambda,\mu)$,
is an arc length parametrisation of a component of $\gamma$. Then $\bt(s)=\gamma'(s)$ is a unit tangent vector and $\bt'(s)=\kappa(s)\bn(s)$,
where $\kappa(s)$ is the Minkowski curvature
of $\gamma$ at $s$ and $\bn$ is the unit normal vector at $s$. The tangent and unit normal vectors are pseudo-orthogonal so they are of
different types, that is, one is spacelike and the other is timelike.  

When $\gamma$ is not necessarily parametrised by arclength, the unit tangent is given by
$$T(t) = \frac{\gamma'(t)}{ |\langle \gamma'(t),\gamma'(t) \rangle |^\frac{1}{2}},$$   
the unit normal by 
$$N(t)= (-1)^\beta T(t),$$
where $\beta=1$ if $\gamma$ is spacelike and $\beta=2$ if $\gamma$ is timelike, 
and the Minkowski curvature (dropping the parameter $t$) is given by $$\kappa = \frac{\langle  \gamma', \gamma''^\perp \rangle }{|\langle \gamma', \gamma' \rangle |^\frac{3}{2}}.$$

\section{The Minkowski Symmetry Set}\label{sec:struc}
The evolute of a spacelike or timelike component of $\gamma$ is the image of the map
$$
e(t)=\gamma(t)-\frac{1}{\kappa(t)}N(t).
$$

In general, the curvature tends to infinity as $t$ tends to $\lambda$ or $\mu$ and the evolute of the curve $\gamma$ is not defined at the
lightlike points. However, the caustic of $\gamma$ is defined everywhere and contains the evolute of $\gamma$ (see for example \cite{ST}). 
The caustic can be defined via the {\it the family of distance-squared functions} $f:S^1 \times \mathbb R^2_1 \to  \mathbb R$
on $\gamma$ given by
$$
f(t,c)=\langle \gamma(t)-c,\gamma(t)-c\rangle.
$$

Denote by $f_c:S^1\to \mathbb R$ the function given by $f_c(t)=f(t,c)$.
We say that $f_c$ has an $A_k$-singularity at $t_0$ if $f_c'(t_0)=f_c''(t_0)=\ldots =f_c^{(k)}(t_0)=0$ and $f_c^{(k+1)}(t_0)\ne 0$.
This is equivalent to
the existence of a local re-parametrisation $h$ of $\gamma$ at $t_0$ such that $(f\circ h)(t)=\pm t^{k+1}$.
Geometrically, $f_c$ has an $A_k$-singularity if and only if the curve $\gamma$ has contact of order $k+1$
at $\gamma(t_0)$ with the pseudo-circle of centre $c$ and radius $r$, with $r=\langle \gamma(t_0)-c, \gamma(t_0)-c \rangle$.
Thus, the curve $\gamma$ has order of contact 1 with a pseudo-circle
at $t_0$ if it intersects transversally the pseudo-circle at $\gamma(t_0)$.
The order of contact is 2 if the circle and the curve have ordinary tangency at $\gamma(t_0)$.

The caustic of $\gamma$ is the local component $\mathcal B_1$ of the bifurcation set of the family $f$, given by
$$
\mathcal B_1=\{c\in  \mathbb R^2_1\,|\, \exists t \in S^1\mbox{ such that }
 f_c'(t)=f_c''(t)=0\}.
$$

This is the set of points $c\in \mathbb R^2_1$ such that the germ
$f_c$ has a degenerate singularity at some point $t$. In \cite{ST} it was shown that the caustic of $\gamma$ is defined at all points
on $\gamma$ including its lightlike points where it is a smooth curve and has ordinary tangency with $\gamma$. 



The multi-local component of the bifurcation set of the family $f$ is defined as
$$
\mathcal B_2=\{c\in  \mathbb R^2_1\,|\, \exists t_1,t_2 \mbox{ such that } t_1\ne t_2, \
 f_c(t_1)=f_c(t_2), \, f_c'(t_1)=f_c'(t_2)=0\}. 
$$
The full-bifurcation set of $f$ is defined as
$$\mbox{\rm Bif}(f)=\mathcal B_1 \cup \mathcal B_2.$$

\begin{defn}
The {\it Minkowski Symmetry Set} (MSS) of $\gamma$ is the locus of centres of pseudo-circles
which are tangent to $\gamma$ in at least two distinct points $p$ and $q$. The pairs of points $p,q$ are called bi-tangent pairs.
\end{defn}

The $MSS$ is precisely the multi-local component $\mathcal B_2$ of the bifurcation set of
the family of distance-squared function $f$ on $\gamma$.

In \cite{RT} it is shown that the singularities which can occur on the MSS for a generic plane curve are $A_1$, $A_2, A_1^3, A_2 A_1$ and $A_3$, and that they are all versally unfolded.  It follows that these singularities are also versally unfolded for a 1-parameter of plane curves.  It can happen for a generic 1-parameter family of plane curves that at isolated points one of the above singularities occurs at lightlike points and this case is also dealt with in \cite{RT}.  It only remains now to show the versality and the transition type for the other generically occurring singularities for a 1-parameter family of plane curves, namely $A_1^4, A_2A_1^2, A_3A_1, A_2^2$ and $A_4$.  

In \cite{BG86} it was shown that for general functions some of  these singularities occur in two distinct transition types.  For example, in the $A_1^4$ case   there exist two types referred to as  $A_1^4(a)$ and $A_1^4(b)$.  It was shown in that paper that only types $A_1^4(a), A_2^2(a), A_2^2(b), A_1A_3(a), A_1^2A_2(a)$ and $A_4$ could occur for (Euclidean) symmetry sets (see Table on \pageref{Table}).  In the present paper a similar analysis is carried out for the Minkowski symmetry set and the geometric conditions for the possible types are determined.  In particular,  the following theorem is proven:
      
\begin{theo}
The possible transitions types of the Minkowski Symmetry set for a generic curve are $A_1^4(a), A_1^4(b), A_2^2(a), A_2^2(b)$, $A_1A_3(a)$, $A_1A_3(b)$, $A_1^2A_2(a)$, $A_1^2A_2(b)$ and $A_4$.
\end{theo}


Each generically occurring singularity type 
 is considered in turn.  Considering the reduction of the distance-squared family to its normal form, the necessary geometrical criteria for each transition type (eg. $a$ or $b$) is determined.

\begin{figure}[tp]
\includegraphics[width=11cm, height=20cm]{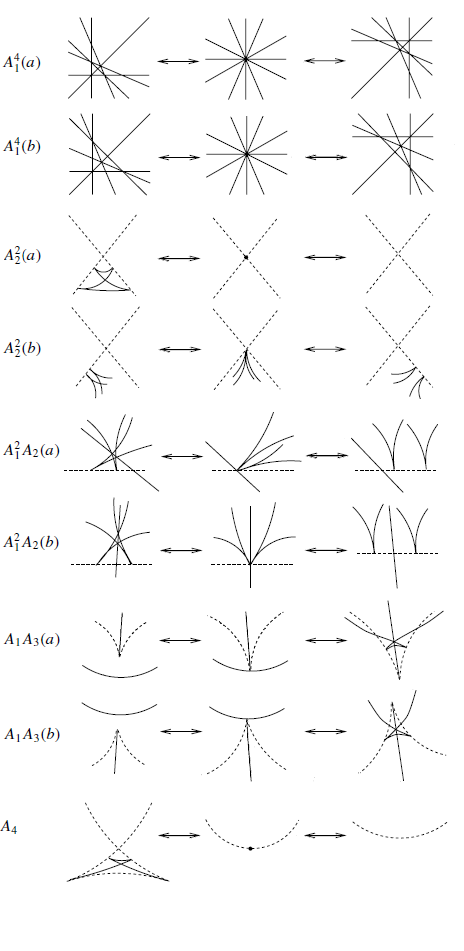}
\caption{The Transitions that can occur on Minkowski Symmetry Sets.}
\label{fig:A14}
\end{figure}






\section{The $A_1^4$ singularity}
Consider the standard multi-versal unfolding of an $A_1^4$ singularity given by
$$G:\mathbb{R}^{(4)} \times \mathbb{R}^3 \to \mathbb{R},$$
where $\mathbb{R}^{(4)}$ denotes the set of parameters $t_1, t_2, t_3, t_4$, $\mathbb{R}^3$ denotes the $(y_1, y_2, y_3)$-space of unfolding parameters and the multi-versal unfolding $G$ is given by 
\begin{eqnarray*}
G_i&:&(t_i, \by) \mapsto t_i^2+ y_i, i=1,2 \ {\rm and} \ 3   \\
G_4&:&(t_4, \by) \mapsto  t_4^2.  
\end{eqnarray*}

Consider now four families of curve segments $\gamma_1, \gm_2, \gm_3$ and $\gm_4$ each being close to one of the tangency points.  With family parameter $u$, denote these segments as $\gm_{i, u} (s_i) = (X_{i, u}(s_i), Y_{i, u}(s_i))$, where the arclength parameters $s_i$ are close to zero.  Take ${\bx} = (x_1, x_2) \in \mathbb R^2_1$, and denote by $\bx_0$ the $A_1^4$-point on the MSS.  Then the family of Minkowski distance functions on the family of curve segments consists of four germs 
$$F_i: \mathbb R \times  \mathbb R \times  \mathbb R^2_1, (0, 0, \bx_0) \to  \mathbb R, $$ 
given by
\begin{eqnarray}
\nonumber
F_i(s_i, u, \bx) = \langle \bx - \gm_{i, u}, \bx - \gm_{i, u} \rangle. 
\label{dis_sqrd_func}
\end{eqnarray}

Using standard techniques, as outlined in \cite{BG86}, and used for example in \cite{HoltTh} and \cite{Pollitt}, is to reduce the family $F_i$ to the standard family $G_i$.  The big bifurcation set (BBS), which sits in $\by$-space and comprises of subsets which correspond to $A_1^2$ sets of $G$, contains all the possible types bifurcations of $A_1^4$, and the individual bifurcation sets can be recovered locally by slicing the BBS with non-singular families of surfaces passing through the origin in $\by$-space.  Firstly, the possible generic transition types and their criteria are found, and then through keeping track of the geometric properties in reducing the family to the standard type, the relevant bifurcation type can be determined.

\subsection{Bad planes}
Following \cite{BG86}, a plane containing the origin given by the equation $a_1 y_1 + a_2 y_2 + a_3 y_3 = 0$
is called a {\it bad plane} if it contains any of the limiting tangent vectors to the strata of the big bifurcation set of $G$.  
Non-generic transitions occur when these slicing surfaces are themselves tangent to the limiting tangent vectors to the strata of the big bifurcation set tending to the origin.      
A plane  can be represented by a point with homogenous coordinates $(a_1 : a_2 : a_3)$ in the real projective plane $\mathbb{R} P^2$ and the pencils of bad planes therefore correspond to lines in $\mathbb{R} P^2$. 
    
If $\Delta$ represents the se of bad planes each component of $\mathbb{R} P^2 - \Delta$ represent collections of normals, which as kernels of $dh(0)$ give $C^0$-stratified equivalent functions of $h$. (For remarks on stratified equivalence see for example \cite{BG86} and \cite{Bru86}.)  Each connects component of $\mathbb{R} P^2 - \Delta$ can potentially give a different type of transition.  By considering each region in turn and identifying the type of  transition it is possible to determine the criteria for realising each one.

The one-dimensional strata adjacent to the BBS for the standard $A_1^4$ are
\begin{eqnarray*}
  A_1^3 &:&  \{ (a_1, a_2, a_3) = (t_1,  t_1, t_1) \cup (t_1, 0, 0)  \cup (0, 0, t_2) \cup (0, 0, t_3)  \}  \\
	A_1^2 / A_1^2 &:&   \{ (a_1, a_2, a_3) = (t_1, t_1, 0) \cup (0, t_2, t_2) \cup (t_3, 0, t_3) \}. 
\end{eqnarray*}

The limiting tangent vectors to these one-dimensional strata are therefore given by 
$(1, 0, 0), (0, 1, 0), (0, 0, 1)$, $(1, 1, 1)$, $(1, 1, 0), (0, 1, 1)$, and $(1, 0, 1)$ so the bad planes are given by
$a_1 = 0$, $a_2 =0$, $a_3 =0$, and $a_1 + a_2 + a_3= 0$, $a_1 + a_2= 0$, $a_2 + a_3= 0$, and $a_1 + a_3= 0$.  

\begin{figure}
\centering
\begin{minipage}{.5\textwidth}
  \centering
  \includegraphics[width=.6\linewidth]{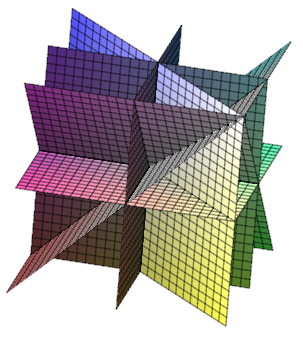}
\end{minipage}%
\begin{minipage}{.5\textwidth}
  \centering
  \includegraphics[width=0.87\linewidth]{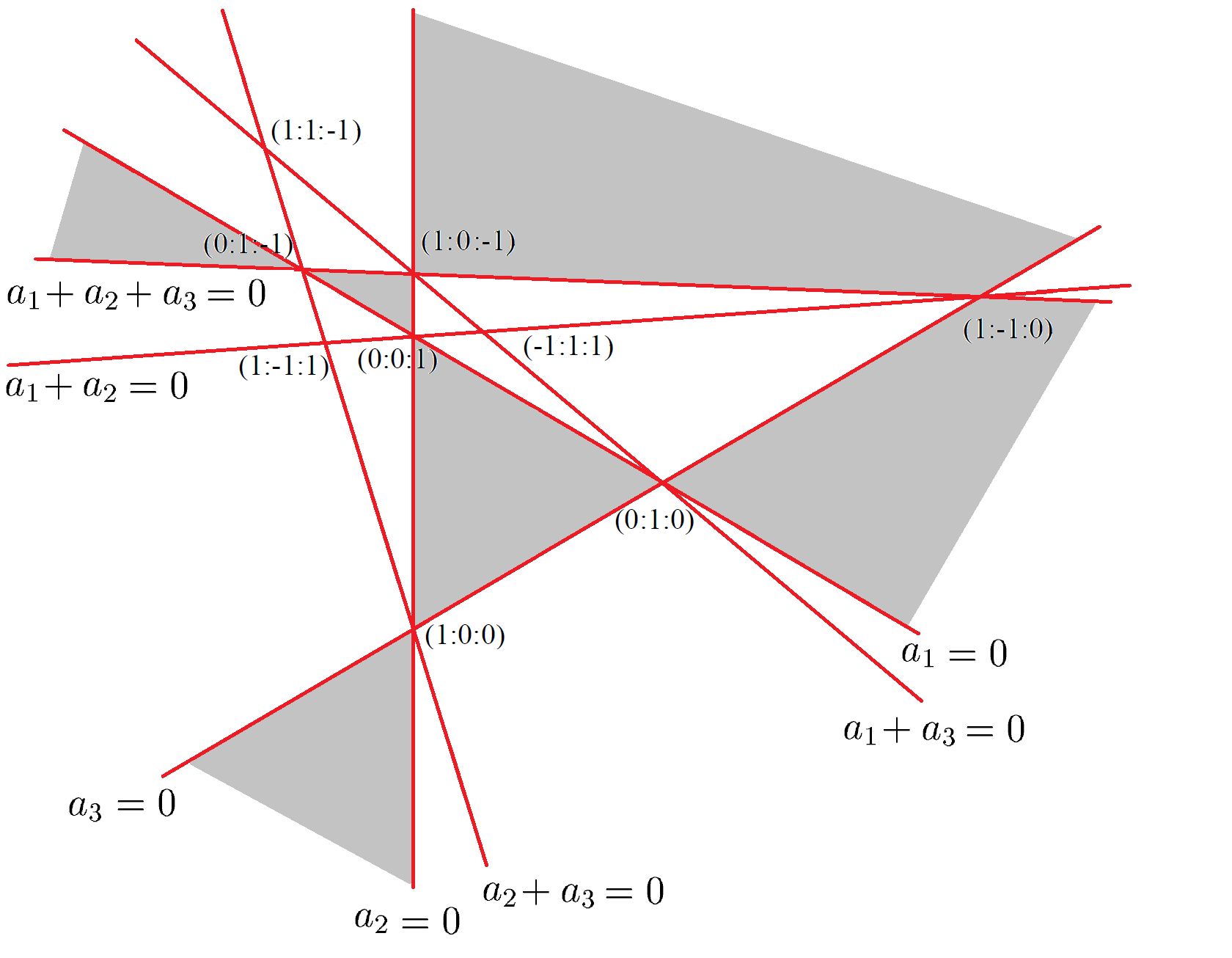}
\end{minipage}
  \captionof{figure}{Left: The BBS for $A_1^4$. Right: The regions determining the different types for $A_1^4$.}
  \label{fig:A14-BBS}
\end{figure}

It is determined that the shaded regions of Figure \ref{fig:A14-BBS} (right) correspond to one type of transition and the non-shaded regions give another from which the following proposition can be deduced.

\begin{prop} 
If $a_1a_2a_3(a_1 + a_2 + a_3)$ is negative the point $(a_1 : a_2 : a_3)$ lies in the shaded region of Figure \ref{fig:A14-BBS} (right) and the corresponding full bifurcation set has type $A_1A_3(a)$.  If however $a_1a_2a_3(a_1 + a_2 + a_3)$ is positive, then the point lies in the unshaded region and the corresponding full bifurcation set is of type $A_1A_3(b)$. 
\end{prop}

Since it is assumed that each $F_i$ is a multi-versal unfolding, then by the uniqueness of multi-versal unfoldings each of the unfoldings $G_i$ in the standard multi-versal unfolding $G$ can be induced from the affine distance functions $F_i$ by
\begin{eqnarray}
G_i(t_i, \by) = F_i(A_i(t_i, \by), B(\by)) + C(\by), {\rm for \ } i= 1, 2, 3 {\rm \ and \ } 4,
\label{diffeo_func}  
\end{eqnarray}
where each $A_i: \mathbb R \times \R^3 \to \R $ is a germ at $(0, {\bf 0})$ and $B, C$ denote the germs $B: (\R^3, {\bf 0}) \to (\R \times \R^2, (0, \bx_0))$  and $C: (\R^3, {\bf 0}) \to (\R, d_0).$ 

$$\begin{CD}
\R \times \R^3   @>G>> \R \times\R^3 @>>> \R^3 @>h>>  \R\\
@VV(A_i \times B)V @VV(-C \times B)V @VVBV  @VV{\rm identity}V \\
C @>F>> D @>>> D @>\pi_1>> \R
\end{CD}$$

\noindent
From the commutative diagram it can be seen that $h =\pi_1 \circ B$, 
where $\pi_1$ denotes projection onto the first coordinate.  Thus, $B_1$ (where $B_i$ denotes the $i^{th}$ component of $B$) is the map $h$ on the standard $A_1^4$ set (the BBS), which corresponds to the plane through the origin in $\by$-space representing the tangent plane to the surface with which we are slicing the BBS.  This tangent plane thus corresponds to the kernel of the map $h$ on the BBS, i.e.
$${\rm ker} \ dB_1: \R^3 \to \R, \ {\rm with \ matrix \ } 
\left. \left( \frac{\partial B_1}{\partial y_1},  \frac{\partial B_2}{\partial y_2},  \frac{\partial B_3}{\partial y_3} \right) \right|_{\by={\bf 0}}. $$
Hence the kernel plane has equation
$$\left. \left. \left. \frac{\partial B_1}{\partial y_1} \right|_{\by={\bf 0}} y_1 +  \frac{\partial B_2}{\partial y_2} \right|_{\by={\bf 0}} y_2 +  \frac{\partial B_3}{\partial y_3} y_3 \right|_{\by={\bf 0}}=0.$$

\begin{prop}
The MSS has a transition of type $A_1^4(a)$ if there are odd number of points on each branch and is of type $A_1^4(b)$ if there are an even number of points on each branch.
\end{prop}

\noindent
{\bf Proof.}
Consider the case $i=1$:
$$\left. \left( \frac{\partial G_1}{\partial t_1}  \ \frac{\partial G_1}{\partial y_1} \ \frac{\partial G_1}{\partial y_2} \ \frac{\partial G_1}{\partial y_3}  \right)   \right|_{\by={\bf 0}}   = (2 t_1 \ 1 \ 0 \ 0). $$

Using relation (\ref{diffeo_func}) and applying the chain rule for derivatives gives:
$$= \left. \left(\frac{\partial F_1}{\partial s_1}  \ \frac{\partial F_1}{\partial x_1} \ \frac{\partial F_1}{\partial x_2} \ \frac{\partial F_1}{\partial x_3} \right) \right|_{(A_1(t_1, {\bf 0}), \bx_0)} \times \left. \left( \begin{array}{cccc} \frac{\partial A_1}{\partial t_1} & \frac{\partial A_1}{\partial y_1} & \frac{\partial A_1}{\partial y_2} & \frac{\partial A_1}{\partial y_3} \\ 
0 & \frac{\partial B_1}{\partial y_1}& \frac{\partial B_1}{\partial y_2} & \frac{\partial B_1}{\partial y_3} \\ 0
& \frac{\partial B_2}{\partial y_1}& \frac{\partial B_2}{\partial y_2} & \frac{\partial B_2}{\partial y_3} \\
0 & \frac{\partial B_3}{\partial y_1}& \frac{\partial B_3}{\partial y_2} & \frac{\partial B_3}{\partial y_3}  \end{array} \right) \right|_{(t_1,{\bf 0})} + \left. \left(0 \ \frac{\partial C}{\partial y_1} \ \frac{\partial C}{\partial y_2} \ \frac{\partial C}{\partial y_3}  \right) \right|_{\by={\bf 0}} $$

The same can be done for $G_2, G_3$ and $G_4$, which have the right side of the first line as $(2t_2 \ 0 \ 1 \ 0)$, $(2t_3 \ 0 \ 0 \ 1)$ and $(2t_4 \ 0 \ 0 \ 0)$ respectively.  Now
$\frac{\partial F_i}{\partial s_i} (0, \bx_0) \equiv 0$ because $F_i$ has an $A_1$ singularity at $(0, \bx_0)$.  Also, $\frac{\partial F_i}{\partial x_1} = - 2x_1 + 2 X_{u, i}(s_i), \frac{\partial F_i}{\partial x_2} = 2x_2 - 2Y_{u, i}(s_i).$   The substitution $t_i=0$ can be made since only the 0-jets are required.

Taking all the $G_i$ together gives the system:
\begin{eqnarray}
\left( \begin{array}{ccc}   1 & 0 & 0 \\ 0 & 1 & 0 \\ 0 & 0 & 1 \\ 0 & 0 & 0
\end{array} \right) = \left. \left( \begin{array}{ccc}   \frac{\partial F_1}{\partial u} & \frac{\partial F_1}{\partial x_1} & \frac{\partial F_1}{\partial x_2} \\ \frac{\partial F_2}{\partial u} & \frac{\partial F_2}{\partial x_1} & \frac{\partial F_2}{\partial x_2} \\ \frac{\partial F_3}{\partial u} & \frac{\partial F_3}{\partial x_1} & \frac{\partial F_3}{\partial x_2} \\
\frac{\partial F_4}{\partial u} & \frac{\partial F_4}{\partial x_1} & \frac{\partial F_4}{\partial x_2}
\end{array} \right) \right|_{(A(t_i, {\bf 0}), \bx_0)} \times JB +  \left( \begin{array}{c}   JC \\ JC  \\ JC \\ JC
\end{array} \right)
\label{eq1_Gi}
 \end{eqnarray}
where, for conciseness, $JB$ and $JC$ denote the matrices
$$JB =  \left. \left( \begin{array}{ccc}   \frac{\partial B_1}{\partial y_1} & \frac{\partial B_1}{\partial y_2} & \frac{\partial B_1}{\partial y_3} \\   \frac{\partial B_2}{\partial y_1} & \frac{\partial B_2}{\partial y_2} & \frac{\partial B_2}{\partial y_3} \\
 \frac{\partial B_3}{\partial y_1} & \frac{\partial B_3}{\partial y_2} & \frac{\partial B_3}{\partial y_3}
 \end{array} \right)  \right|_{\by=0}, \ \ \  JC= \left. \left( \begin{array}{ccc}   \frac{\partial C}{\partial y_1} & \frac{\partial C}{\partial y_2} & \frac{\partial C}{\partial y_3} 
 \end{array} \right)\right|_{\by=0}.
$$
Subtracting the bottom row from the other rows in equation (\ref{eq1_Gi}) gives
$$\left( \begin{array}{ccc}   1 & 0 & 0 \\ 0 & 1 & 0 \\ 0 & 0 & 1 \\ 0 & 0 & 0
\end{array} \right) = \left. \left( \begin{array}{ccc}   \frac{\partial F_1}{\partial u} - \frac{\partial F_4}{\partial u}& \frac{\partial F_1}{\partial x_1} - \frac{\partial F_4}{\partial x_1} & \frac{\partial F_1}{\partial x_2} - \frac{\partial F_4}{\partial x_2} \\ \frac{\partial F_2}{\partial u} - \frac{\partial F_4}{\partial u} & \frac{\partial F_2}{\partial x_1} - \frac{\partial F_4}{\partial x_1} & \frac{\partial F_2}{\partial x_2} - \frac{\partial F_4}{\partial x_2} \\ \frac{\partial F_3}{\partial u} - \frac{\partial F_4}{\partial u} & \frac{\partial F_3}{\partial x_1} - \frac{\partial F_4}{\partial x_1}& \frac{\partial F_3}{\partial x_2} - \frac{\partial F_4}{\partial x_2}\\
\frac{\partial F_4}{\partial u} & \frac{\partial F_4}{\partial x_1} & \frac{\partial F_4}{\partial x_2}
\end{array} \right) \right|_{(A(t_i, {\bf 0}), \bx_0)} \times JB +  \left( \begin{array}{c}  {\bf 0 } \\ {\bf 0 }  \\ {\bf 0 } \\ $JC$ \end{array} \right).
 $$
Substituting $\frac{\partial F_i}{\partial x_1}$ and $\frac{\partial F_i}{\partial x_2}$ and ignoring the last row yields the following system:

$$I_3 =  \left( \begin{array}{ccc}   \frac{\partial F_1}{\partial u} - \frac{\partial F_4}{\partial u}& X_1 - X_4 & -Y_1 + Y_4 \\ \frac{\partial F_2}{\partial u} - \frac{\partial F_4}{\partial u}  & X_2 - X_4 &  -Y_2 + Y_4 \\   \frac{\partial F_3}{\partial u} - \frac{\partial F_4}{\partial u}   &  X_3 - X_4  &  -Y_3 + Y_4 \\
 \end{array} \right) \times \left( \begin{array}{ccc}   \frac{\partial B_1}{\partial y_1} & \frac{\partial B_1}{\partial y_2} & \frac{\partial B_1}{\partial y_3} \\   \frac{\partial B_2}{\partial y_1} & \frac{\partial B_2}{\partial y_2} & \frac{\partial B_2}{\partial y_3} \\
 \frac{\partial B_3}{\partial y_1} & \frac{\partial B_3}{\partial y_2} & \frac{\partial B_3}{\partial y_3}
 \end{array} \right)
$$
where $I_3$ represents the $(3 \times 3)$ identity matrix.

The derivatives of $B_1$ can now be evaluated.  Since the product of the two matrices is the identity, they must be inverse to each other.  Now, the inverse of the first matrix can be used to calculate the required entries of the second matrix.  So, 
$$\frac{\partial B_1}{\partial y_1} = \beta {\rm det} \left(\begin{array}{cc}  X_2 - X_4 &  -Y_2 + Y_4  \\    X_3 - X_4  &  -Y_3 + Y_4 \end{array} \right).$$
Multiplying the second column by $-1$ gives

$$\frac{\partial B_1}{\partial y_1} = - \beta {\rm det} \left(\begin{array}{cc}  X_2 - X_4 &  Y_2 - Y_4  \\    X_3 - X_4  &  Y_3 - Y_4 \end{array} \right).$$

Similarly, 
$$\frac{\partial B_1}{\partial y_2} = - \beta {\rm det} \left(\begin{array}{cc}  X_1 - X_4 &  Y_1 - Y_4  \\    X_3 - X_4  &  Y_3 - Y_4 \end{array} \right),$$
$$\frac{\partial B_1}{\partial y_3} = - \beta {\rm det} \left(\begin{array}{cc}  X_1 - X_4 &  Y_1 - Y_4  \\    X_2 - X_4  &  Y_2 - Y_4 \end{array} \right).$$


Let $q_1 = \gm_2 - \gm_3$,  $q_2 = \gm_3 - \gm_4$, $q_3 = \gm_4 - \gm_1$ and $q_4 = \gm_1 - \gm_2$.
Now, 
$\frac{\partial B_1}{\partial y_1} = -\beta {\rm det} \left(\genfrac{}{}{0pt}{}{q_1}{q_2} \right)$,
$\frac{\partial B_1}{\partial y_2} = \beta {\rm det} \left(\genfrac{}{}{0pt}{}{q_2}{q_3} \right)$,
$\frac{\partial B_1}{\partial y_3} = -\beta {\rm det} \left(\genfrac{}{}{0pt}{}{q_3}{q_4} \right)$,
and 
$\frac{\partial B_1}{\partial y_1} + \frac{\partial B_1}{\partial y_2} + \frac{\partial B_1}{\partial y_3} =  - \beta {\rm det} \left(\genfrac{}{}{0pt}{}{q_3}{q_4} \right) $.

\begin{figure}[tp]
\includegraphics[width=13.1cm, height=4.1cm]{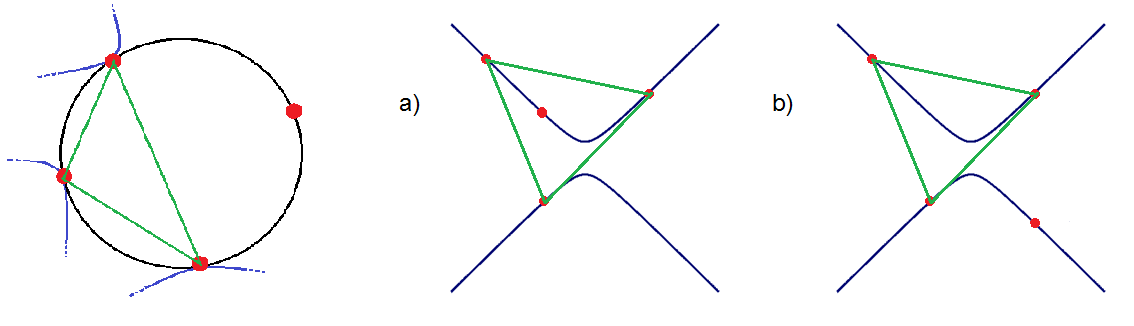}
\caption{Left: Given three points on a circle, a fourth point necessarily lies outsied the triangle formed by the other three.  Right: Given three points on a pseudo-circle, a fourth point can either lie inside (resulting in singularity $A_1^4(a)$), or outside the triangle formed by the other three (resulting in the singularity $A_1^4(b))$.}
\label{fig:A14-conds}
\end{figure}

Now ${\rm det}(q_i, q_j)>0$ if and only if the anticlockwise (Euclidean) angle from $q_i$ to $q_j$ is less than $\pi$.  It then follows that $\frac{\partial B_1}{\partial y_1}  \frac{\partial B_1}{\partial y_2}  \frac{\partial B_1}{\partial y_3} \left( \frac{\partial B_1}{\partial y_1} + \frac{\partial B_1}{\partial y_2} + \frac{\partial B_1}{\partial y_3} \right) >0$ if and only if no point $p_i$ is inside the triangle formed by the other three $p_j$.  This condition fails if and only if there are an even number of points on each branch and the resulting singularity is of type $A_1^4(b)$.  On the other hand, if one of the branches contains only one point, and the other branch contains three, then the triangle formed by the point on the first branch and the `outer' two points of the branch of three will necessarily contain the fourth point (see figure \ref{fig:A14-conds}) and the singularity will be of type $A_1^4(a)$.




\section{The $A_2^2$ singularity}
Consider the following standard multi-versal unfolding of an $A_1^2A_2$ singularity given by
$$G: \mathbb{R}^{(2)} \times \mathbb{R}^3 \to \R $$ 
where $\mathbb{R}^{(2)}$ denotes the parameters $t_1, t_2$ and $\mathbb{R}^3$ denotes the unfolding parameters $\ba = (a_1, a_2, a_3)$ and the multi-versal unfolding is given by the two unfoldings:
\begin{eqnarray*}
G_1(t_1, \ba) &=& t_1^3 + a_1t_1 + a_2, \\ 
G_2(t_2, \ba) &=& t_2^3 + a_3t_2.
\end{eqnarray*}

\subsection{The bad planes}
The one-dimensional strata adjacent to $A_2^2$  are
\begin{eqnarray*}
  A_1 A_2  &:&  \{ (a_1, a_2, a_3) = (-3t_1^2,  2t_1^3, 0) \cup (0, -2 t_2^3, -3t_2^2)   \}  \\
	A_1^2 / A_1^2 &:&   \{ (a_1, a_2, a_3) = (-3t_2^2, 0, -3t_2^2) \}. 
\end{eqnarray*}

The limiting tangent vectors to these one-dimensional strata are given by 
$(1, 0, 0), (0, 1, 0)$  and $(1, 1, 0)$ so the bad planes are given b
$a_1 = 0$, $a_3 =0$ and $a_1 + a_3= 0$.

\begin{figure}
\centering
\begin{minipage}{.5\textwidth}
  \centering
  \includegraphics[width=0.7\linewidth]{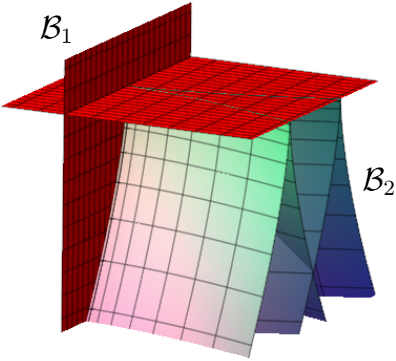}
\end{minipage}%
\begin{minipage}{.5\textwidth}
  \centering
  \includegraphics[width=0.97\linewidth]{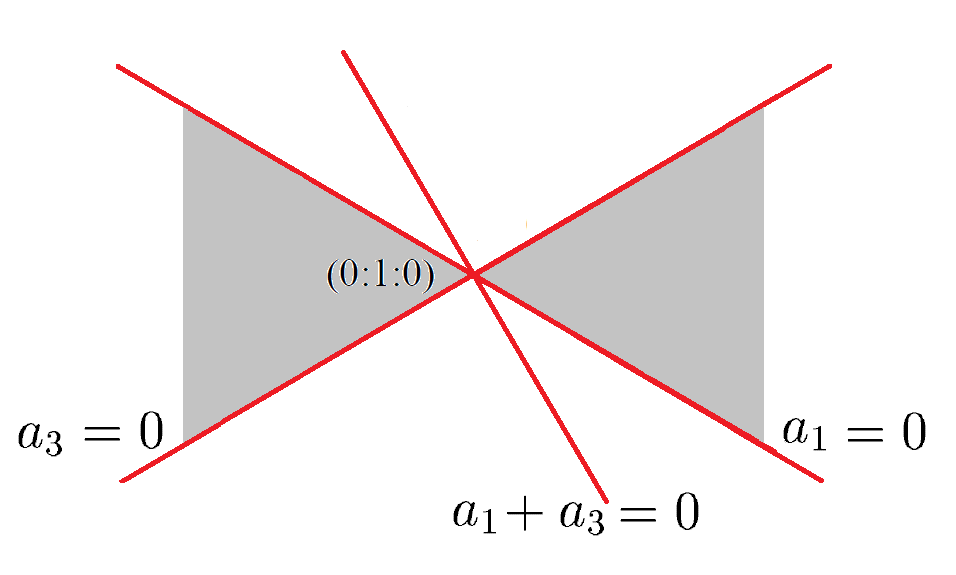}
\end{minipage}
\captionof{figure}{Left: The BBS for $A_2^2$. Right: The regions determining the different types for $A_2^2$.}
\label{fig:A22-BBS}
\end{figure}

Similarly to the previous case, the following proposition can be deduced.
\begin{prop} 
If $a_1a_2$ is negative the point $(a_1 : a_2 : a_3)$ lies in the unshaded region of Figure \ref{fig:A22-BBS} (right) and the corresponding full bifurcation set has type $A_2^2(a)$.  If however $a_1a_3$ is positive, then the point lies in the shaded region and the corresponding full bifurcation set is of type $A_2^2(b)$. 
\end{prop}

The Minkowski distance function on the two curve segments near the $A_2^2$ points consists of the two germs
$$F_1(t_1, u, x) = \langle \gamma_1(t_1, u) - x, \gamma_1(t_1, u) - x \rangle $$ 
$$F_2(t_2, u, x) = \langle \gamma_2(t_2, u) - x, \gamma_2(t_2, u) - x \rangle.$$ 
To reduce to $G_1$ and $G_2$, as in the $A_1^4$ case, using (\ref{diffeo_func}) and applying the 
chain rule gives the system:
$$\left( \begin{array}{ccc}   1 & 0 & 0 \\ 0 & 1 & 0 \\ 0 & 0 & 1 \\ 0 & 0 & 0
\end{array} \right) = \left. \left( \begin{array}{ccc}   \frac{\partial^2 F_1}{\partial t_1 \partial u} & \frac{\partial^2 F_1}{\partial t_1 \partial x_1} & \frac{\partial^2 F_1}{\partial t_1 \partial x_2} \\ 
\frac{\partial F_1}{\partial u} & \frac{\partial F_1}{\partial x_1} & \frac{\partial F_1}{\partial x_2} \\ 
\frac{\partial^2 F_2}{\partial t_2 \partial u} & \frac{ \partial^2 F_2}{\partial t_2 \partial x_1} & \frac{\partial^2 F_2}{\partial t_2 \partial x_2}  \\
\frac{\partial F_2}{\partial u} & \frac{\partial F_2}{\partial x_1} & \frac{\partial F_2}{\partial x_2} 
\end{array} \right) \right|_{(A(t_i, {\bf 0}), \bx_0)} \times JB +  \left( \begin{array}{c}   {\bf 0 }  \\ JC  \\ {\bf 0 }  \\ JC
\end{array} \right).
 $$

Subtracting the bottom row from the second and then ignoring the bottom yields
$$\left( \begin{array}{ccc}   1 & 0 & 0 \\ 0 & 1 & 0 \\ 0 & 0 & 1 
\end{array} \right) = \left. \left( \begin{array}{ccc}   \frac{\partial^2 F_1}{\partial t_1 \partial u} & \frac{\partial^2 F_1}{\partial t_1 \partial x_1} & \frac{\partial^2 F_1}{\partial t_1 \partial x_2} \\ 
\frac{\partial F_1}{\partial u} - \frac{\partial F_2}{\partial u} & \frac{\partial F_1}{\partial x_1} - \frac{\partial F_2}{\partial x_1} & \frac{\partial F_1}{\partial x_2} - \frac{\partial F_2}{\partial x_2}  \\ 
\frac{\partial^2 F_2}{\partial t_2 \partial u} & \frac{ \partial^2 F_2}{\partial t_2 \partial x_1} & \frac{\partial^2 F_2}{\partial t_2 \partial x_2}  
\end{array} \right) \right|_{(A(t_i, {\bf 0}), \bx_0)} \times JB 
 $$

 

We can write
$A_i(t_i, 0)=\alpha_i t_i + {\rm higher \ terms}$  where
$$\alpha_i = (- \kappa / \kappa_i)^{\frac{1}{3}}$$
and here $\kappa$ is the Minkowski curvature of $\gamma$ at the two points of contact and $\kappa_i'$ is the derivative of Minkowski curvature with respect to arclength on $\gamma$.

Differentiating $F_1$ (for example, though the same applies for $F_2$) gives 
$$\frac{1}{2}\frac{\partial F_1(A_1(t_1, u), x)}{\partial t_1 } = \alpha_1\langle (\gamma_1(t_1, u) - x), T_1 \rangle$$
and differentiating this with respect to $x$ gives
$$\left(\frac{1}{2}\frac{\partial^2 F_1(A_1(t_1, u), x)}{\partial t_1 \partial x_1 }, \frac{1}{2}\frac{\partial^2 F_1(A_1(t_1, u), x)}{\partial t_1 \partial x_2 }\right) =  \alpha_1 (X_1', -Y_1').$$  

For the middle row we have
$$\left(\frac{1}{2}\frac{\partial F_i(A_i(t_i, u), x)}{\partial x_1}, \frac{1}{2}\frac{\partial F_i(A_i(t_i, u), x)}{\partial x_2}\right) =  \langle (\gamma_1(t, u) - x),  (-1, -1)\rangle.$$ 
Since $F_1$ has an $A_2$ singularity, $(\gamma(t, u) - x)$ can be written as $\frac{1}{\kappa_M} N_M$ and substituting this yields
$$\left(\frac{1}{2}\frac{\partial F_i(A_i(t_i, u), x)}{\partial x_1 }, \frac{1}{2}\frac{\partial F_i(A_i(t_i, u), x)}{\partial x_2}\right)  =  2 \frac{1}{\kappa_M}(Y'_i, -X'_1).$$

Substituting these derivatives into the matrix equation gives:
$$\left( \begin{array}{ccc}   1 & 0 & 0  \\ 0 & 1 & 0 \\ 0 & 0 & 1
\end{array} \right) = \left. \left( \begin{array}{ccc}  * & 2\alpha_1 X_1' & -2\alpha_1 Y_1'   \\ 
* & \frac{2}{\kappa}(Y'_1 -Y'_2) & \frac{2}{\kappa}(X'_2 -X'_1) \\
* & 2\alpha_2 X_2' & -2\alpha_2 Y_2'  
\end{array} \right) \right|_{(A(t_i, {\bf 0}), \bx_0)} \times JB. 
 $$


Evaluating the cofactors gives $\frac{\partial B_1}{\partial a_1}=  \frac{4}{\kappa} \alpha_2 (\langle T_1, T_2 \rangle \pm 1)$ 
and$\frac{\partial B_1}{\partial a_3} = \frac{4}{\kappa} \alpha_1 (\langle T_1, T_2 \rangle \pm 1)$, where the sign of $\pm$ is the same for both derivatives and depends on whether the curves are spacelike or timelike.

The type of transition that occurs depends on the sign of $\frac{\partial B_1}{\partial a_1}\frac{\partial B_1}{\partial a_3}$
Now $\frac{\partial B_1}{\partial a_1}\frac{\partial B_1}{\partial a_3}=\frac{8}{\kappa_M^2}\alpha_1\alpha_2(\langle T_1, T_2 \rangle \pm 1)^2$ so the sign, and hence the transition type, depends on whether $\kappa'_1 \kappa'_2$ is positive or negative.

\begin{prop}
In the multi-versal $A_2^2$ situation, assume in addition to $\kappa_i' \neq 0$, that $\kappa_1' + \kappa_2' \neq 0$ ($\kappa_i' =$ the derivative of curvature on $\gamma_0$ with respect to arclength at the two contact points).
Then the $A_2^2(a)$ or ``moth transition" occurs when $\kappa_1'\kappa_2' > 0$ and the $A^2_2(b)$ or ``nib transition'' occurs  when  $\kappa_1'\kappa_2' < 0$.
\end{prop}


\section{The $A_1^2A_2$ singularity}
\label{sectionA12A2}
Consider the following standard multi-versal unfolding of an $A_1^2A_2$ singularity given by
$$G: \mathbb{R}^{(3)} \times \mathbb{R}^3 \to \R $$ 
where $\mathbb{R}^{(3)}$ denotes the parameters $t_1, t_2, t_3$ and $\mathbb{R}^3$ denotes the unfolding parameters $\ba = (a_1, a_2, a_3)$ and the multi-versal unfolding is given by the two unfoldings:
\begin{eqnarray*}
G_1(t_1, \ba) &=& t_1^3 + a_1 t_1, \\ 
G_2(t_2, \ba) &=& t_2^2 + a_2, \\
G_3(t_3, \ba) &=& t_3^2 + a_3. 
\end{eqnarray*}

\subsection{The big bifurcation set}
At an $A_1^2A_2$ point the $\mathcal B_2$ set consists of three parts:
The first is given as the solution of $G_1 = G_2$ and $G_1' = G_2' = 0$ and is a semi-cubic cylinder with the parametrisation $\{[-3t_1^2,  2t_1^3, a_3]\}$. The second is given as the solution of $G_1 = G_3$ and $G_1' = G_3' = 0$ and is a semi-cubic cylinder with the parametrisation $(-3 t_1^2,  a_2, 2t_1^3)$.  The third component is a smooth surface which is the solution set of $G_2= G_3$ and   $G_2' = G_3' = 0$ and can be parametrized as $(a_1, a_2, a_2)$.  The $\mathcal B_1$ component given by $G_1' = G_1'' = 0$ is the smooth surface $(0, a_2, a_3)$.  See Figure \ref{fig:A1A22-BBS} (Left).

\subsection{The bad planes}
The one-dimensional strata adjacent to $A_1^2A_2$  are
\begin{eqnarray*}
  A_1 A_2  &:&  \{ (a_1, a_2, a_3) = (0,  a_2, 0) \cup (0, 0, a_3)   \}  \\
	A_1^3 &:&    \{ (a_1, a_2, a_3) = (-3 t_1^2, -2 t_1^3, -2 t_1^3) \} \\
	A_1^2 / A_1^2 &:&   \{ (a_1, a_2, a_3) = (3t_1^2, 2 t_1^3, -2t_1^3) \} 
\end{eqnarray*}

The limiting tangent vectors to these one-dimensional strata are given by 
$(0, 1, 0)$, $(0, 0, 1)$ and $(1, 0, 0)$ so the bad planes are given by
$a_2 = 0$,  $a_3 =0$ and $a_1 = 0$.

\begin{figure}
\centering
\begin{minipage}{.5\textwidth}
  \centering
  \includegraphics[width=0.64\linewidth]{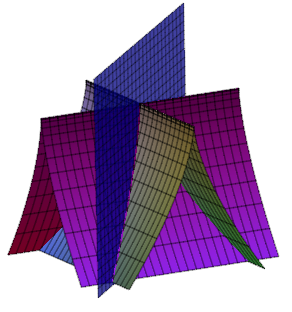}
\end{minipage}%
\begin{minipage}{.5\textwidth}
  \centering
  \includegraphics[width=0.64\linewidth]{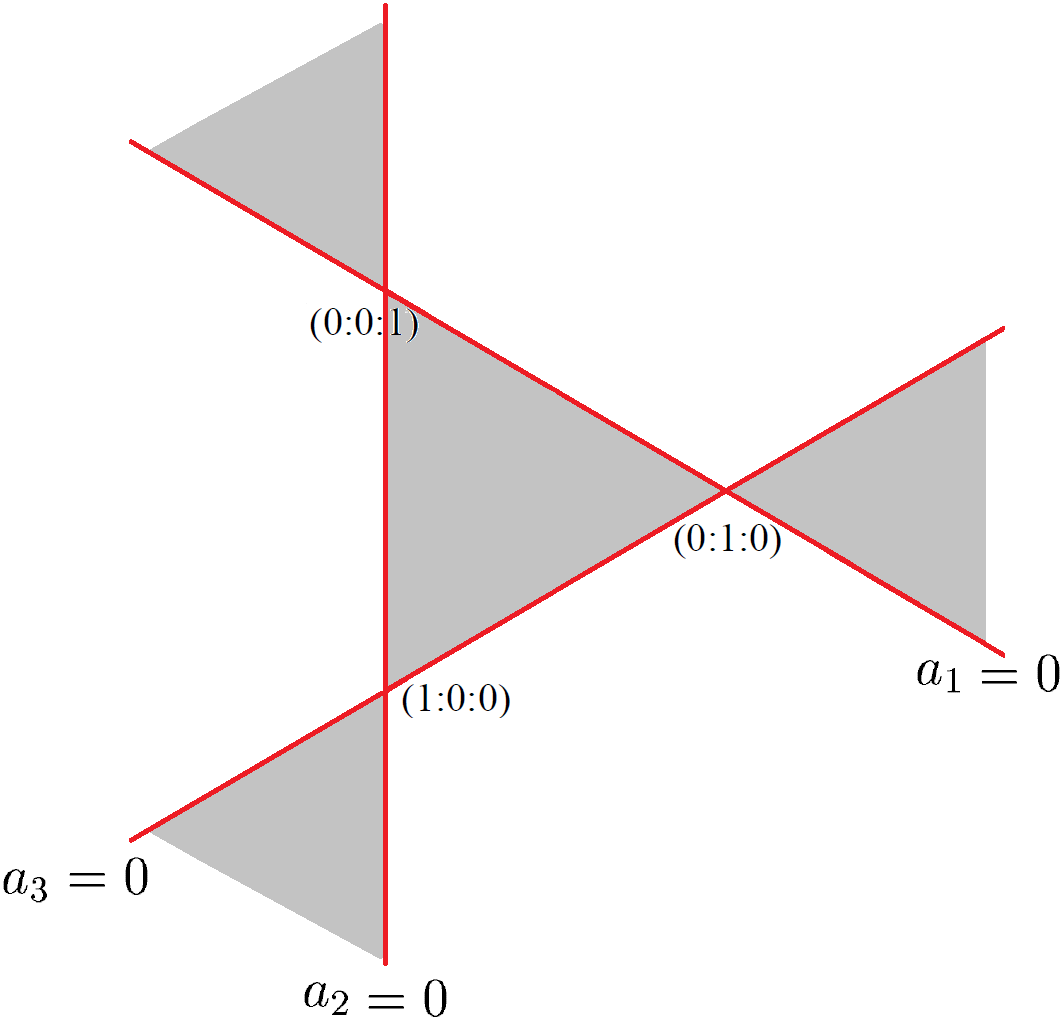}
\end{minipage}
\captionof{figure}{Left: The set $\mathcal B_2$ for $A_1^2A_2$.  The $\mathcal B_1$ set (not shown) is the  plane that contains both cuspidal edges of $\mathcal B_2$. Right: The regions determining the different types for $A_1^2A_2$.}
\label{fig:A1A22-BBS}
\end{figure}

\begin{prop}
If $a_1a_3$ is positive the point $(a_1 : a_2 : a_3)$ lies in the shaded region of Figure \ref{fig:A1A22-BBS} (right) and the corresponding full bifurcation set has type $A_1A_3(a)$.  If however $a_1a_3$ is negative, then the point lies in the unshaded region and the corresponding full bifurcation set is of type $A_1A_3(b)$. 
\end{prop}

Applying the chain rule to (\ref{diffeo_func}) in this case gives the system:
$$\left( \begin{array}{ccc}   1 & 0 & 0 \\ 0 & 0 & 0 \\ 0 & 1 & 0 \\ 0 & 0 & 1
\end{array} \right) = \left. \left( \begin{array}{ccc}   \frac{\partial^2 F_1}{\partial t_1 \partial u} & \frac{\partial^2 F_1}{\partial t_1 \partial x_1} & \frac{\partial^2 F_1}{\partial t_1 \partial x_2} \\ 
\frac{\partial F_1}{\partial u} & \frac{\partial F_1}{\partial x_1} & \frac{\partial F_1}{\partial x_2} \\ 
\frac{\partial F_2}{\partial u} & \frac{ \partial F_2}{ \partial x_1} & \frac{\partial F_2}{ \partial x_2}  \\
\frac{\partial F_3}{\partial u} & \frac{\partial F_3}{\partial x_1} & \frac{\partial F_3}{\partial x_2} 
\end{array} \right) \right|_{(A(t_i, {\bf 0}), \bx_0)} \times JB +  \left( \begin{array}{ccc}  {\bf 0 }  \\ JC \\ JC \\ JC
\end{array} \right).
 $$

Subtracting the second row from the third and fourth rows gives:
$$\left( \begin{array}{ccc}   1 & 0 & 0 \\ 0 & 0 & 0 \\ 0 & 1 & 0 \\ 0 & 0 & 1
\end{array} \right) = \left. \left( \begin{array}{ccc}   \frac{\partial^2 F_1}{\partial t_1 \partial u} & \frac{\partial^2 F_1}{\partial t_1 \partial v_1} & \frac{\partial^2 F_1}{\partial t_1 \partial v_2} \\ 
\frac{\partial F_1}{\partial u}  & \frac{\partial F_1}{\partial x_1}  & \frac{\partial F_1}{\partial x_2} \\ 
\frac{\partial F_2}{\partial u} - \frac{\partial F_1}{\partial u} & \frac{ \partial F_2}{ \partial x_1} - \frac{\partial F_1}{\partial v_1} & \frac{\partial F_2}{ \partial x_2} - \frac{\partial F_1}{\partial x_2} \\
\frac{\partial F_3}{\partial u} - \frac{\partial F_3}{\partial u} & \frac{\partial F_3}{\partial x_1} - \frac{\partial F_1}{\partial x_1} & \frac{\partial F_3}{\partial x_2} - \frac{\partial F_1}{\partial x_2} 
\end{array} \right) \right|_{(A(t_i, {\bf 0}), \bx_0)} \times JB +  \left( \begin{array}{ccc}  {\bf 0 }  \\ JC \\ {\bf 0 } \\ {\bf 0 }
\end{array} \right).
 $$

Ignoring the second row and substituting the derivatives gives
$$\left( \begin{array}{ccc}   1 & 0 & 0  \\ 0 & 1 & 0 \\ 0 & 0 & 1
\end{array} \right) = \left. \left( \begin{array}{ccc} 
* & 2\alpha_1 X_1' & -2\alpha_1 Y_1'  \\ 
* & \frac{2}{\kappa}(Y'_2 -Y'_1) & \frac{2}{\kappa}(X'_1 -X'_2) \\
* & \frac{2}{\kappa}(Y'_3 -Y'_1) & \frac{2}{\kappa}(X'_1 -X'_3)
\end{array} \right) \right|_{(A(t_i, {\bf 0}), \bx_0)} \times JB. 
 $$

Since the bifurcation type depends on whether $\frac{\partial B_1}{\partial a_2}\frac{\partial B_1}{\partial a_3}$ is positive or negative, evaluating these terms using the cofactors of the matrix gives 
\begin{eqnarray}
\nonumber
\frac{\partial B_1}{\partial a_2}{\frac{\partial B_1}{\partial a_3}} &=& 
\frac{16 \alpha_1^2}{\kappa^2}(X_1'^2 - Y_1'^2 -X_1'X_2' + Y_1' Y_2')(X_1'X_3' - X_1'^2 + Y_1'^2 - Y_1' Y_3') \\ 
\nonumber
&=& -\frac{16 \alpha_1^2}{\kappa^2} (\langle T_1, T_1 \rangle - \langle T_1, T_2 \rangle) (\langle T_1, T_1 \rangle - \langle T_1, T_3 \rangle). \\
\label{A12A2condn}
&=& -\frac{16 \alpha_1^2}{\kappa^2} ((-1)^{\beta+1} - \langle T_1, T_2 \rangle) ((-1)^{\beta+1} - \langle T_1, T_3 \rangle).
\end{eqnarray}

If the curves corresponding to the $A_1^2A_2$ point are spacelike, then the pseudo-circle is of type $S^1_1(c, r)$ (radius $r$ and centred at $c$) and can be parametrised as $S^1_1(\theta) = c + r(\cosh(\theta), \pm \sinh(\theta))$, where the $\pm$ allows for the covering of both branches. The unit tangent vectors at $\gamma_i$ are then given by $T_i = (\sinh(\theta_i), \pm \cosh(\theta_i))$.   
If both $\gamma_1$ and $\gamma_i$ ($i=2$ or $3$) lie on the same branch, then 
$$\langle T_1, T_i \rangle = -\sinh(\theta_1)\sinh(\theta_i) + \cosh(\theta_1)\cosh(\theta_i) = \cosh(\theta_1 - \theta_i)$$
so is greater than 1.
If however $\gamma_1$ and $\gamma_i$ lie on opposite branches then 
$$\langle T_1, T_i \rangle = -\sinh(\theta_1)\sinh(\theta_i) - \cosh(\theta_1)\cosh(\theta_i) = -\cosh(\theta_1 + \theta_i)$$
so is less than $-1$.  Since the curves $\gamma_i$ are locally spacelike, $\beta=1$ and the expression (\ref{A12A2condn}) is positive if $\gamma_2$ and $\gamma_3$, that is the two $A_1$ points, lie on the same branch and negative if they lie on opposite branches.  It can be shown that the same result holds if the points are timelike.  It follows that the point is of type $A_1^2A_2$ if of type $(a)$ if the two $A_1$ points lie on the same branch, and of type $(b)$ if they lie on opposite branches of the pseudo-circle.

\section{The $A_3A_1$ singularity}
Consider the following standard multi-versal unfolding of an $A_3A_1$ singularity given by
$$G: \mathbb{R}^{(2)} \times \mathbb{R}^3 \to \R $$ 
where $\mathbb{R}^{(2)}$ denotes the parameters $t_1, t_2$ and $\mathbb{R}^3$ denotes the unfolding parameters $\ba = (a_1, a_2, a_3)$ and the multi-versal unfolding is given by the two unfoldings:
\begin{eqnarray*}
G_1(t_1, \ba) &=& t_1^4 + a_1 t_1^2 + a_2 t_1 + a_3, \\
G_2(t_2, \ba) &=& t_2^2.
\end{eqnarray*}

\subsection{The big bifurcation set}
At an $A_3A_1$ point the $\mathcal B_2$ set itself consists of two parts:
The first is given as the solution to both $G_1 = G_2$ and $G_1' = G_2' = 0$ and is the swallowtail surface parametrised by 
$(a_1, - 4 t_1^3 - 2a_1 t_1, 3 t_1^4 + 2t_1^2a_1)$.  The second component occurs locally near the $A_3$ point and is given by $G_1(t_1) = G_1(-t_1)$ and   $G_1(t_1)' = G_1(-t_1)' = 0$.  This second component is the half plane $(-2t_1^2, 0, y_3)$.  The $\mathcal B_1$ component given by $G_1' = G_1''= 0$ is the semi-cubic cylinder $(- 6 t_1^2, 8 t_1^3, a_3)$, (see Figure \ref{fig:A3A1-BBS} (left)).

\subsection{The bad planes}
The adjacent singularities of codimension 1 are as follows:

\begin{eqnarray*}
  A_3   &:&  \{ (a_1, a_2, a_3) = (0, 0, a_3)  \}    \\
	A_2A_1   &:& \{ (a_1, a_2, a_3) = (-6 t_1^2, 8 t_1^3, -3 t_1^4)  \} \\
	A_1^3 &:& \{ (a_1, a_2, a_3) = 		(-2 t_1^2, 0, t_1^4) \} \\
	A_1^2 / A_1^2   &:& \{ (a_1, a_2, a_3)   = (a_1, 0, 0)  \}
\end{eqnarray*}

The limiting tangent vectors to these one-dimensional strata are given by 
$(1, 0, 0)$, and $(0, 0, 1)$ so the bad planes are given by
$a_1 = 0$ and $a_3 =0$.

\begin{figure}
\centering
\begin{minipage}{.5\textwidth}
  \centering
  \includegraphics[width=0.67\linewidth]{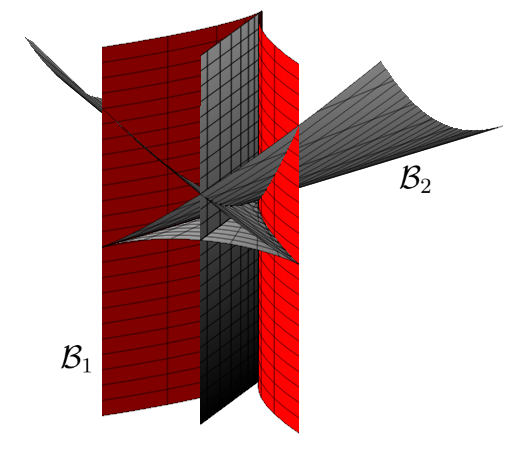}
\end{minipage}%
\begin{minipage}{.5\textwidth}
  \centering
  \includegraphics[width=0.91\linewidth]{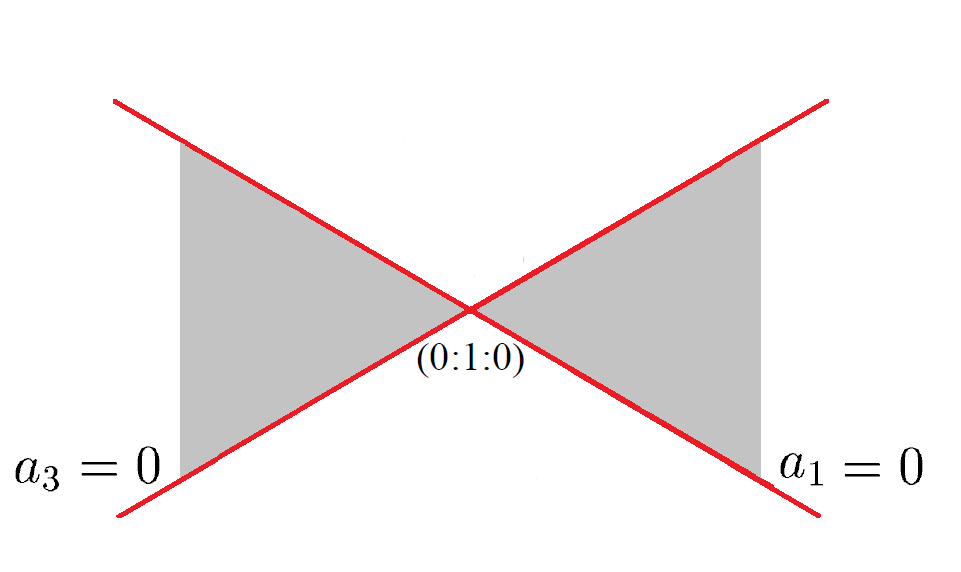}
\end{minipage}
\captionof{figure}{Left: The BBS for $A_3A_1$. Right: The regions determining the different types for $A_3A_1$.}
\label{fig:A3A1-BBS}
\end{figure}

\begin{prop}
If $a_1a_3$ is positive the point $(a_1 : a_2 : a_3)$ lies in the shaded region of Figure \ref{fig:A3A1-BBS} (right) and the corresponding full bifurcation set has type $A_1A_3(a)$.  If however $a_1a_3$ is negative, then the point lies in the unshaded region and the corresponding full bifurcation set is of type $A_1A_3(b)$. 
\end{prop}

Applying the chain rule  to \ref{diffeo_func} gives the system:
$$\left( \begin{array}{ccc}   1 & 0 & 0 \\ 0 & 1 & 0 \\ 0 & 0 & 1 \\ 0 & 0 & 0
\end{array} \right) = \left. \left( \begin{array}{ccc}   
\frac{\partial^3 F_1}{\partial^2 t_1 \partial u} & \frac{\partial^3 F_1}{\partial^2 t_1 \partial x_1} & \frac{\partial^3 F_1}{\partial^2 t_1 \partial x_2}  \\
\frac{\partial^2 F_1}{\partial t_1 \partial u} & \frac{\partial^2 F_1}{\partial t_1 \partial x_1} & \frac{\partial^2 F_1}{\partial t_1 \partial x_2} \\ 
\frac{\partial F_1}{\partial u} & \frac{\partial F_1}{\partial x_1} & \frac{\partial F_1}{\partial x_2} \\ 
\frac{\partial F_2}{\partial u} & \frac{ \partial F_2}{ \partial x_1} & \frac{\partial F_2}{ \partial x_2}  
\end{array} \right) \right|_{(A(t_i, {\bf 0}), \bx_0)} \times JB +  \left( \begin{array}{ccc}  {\bf 0 }  \\ {\bf 0 }\\ JC \\ JC
\end{array} \right).
 $$
Subtracting the last row from the third, and then ignoring the last gives:
$$\left( \begin{array}{ccc}   1 & 0 & 0 \\ 0 & 1 & 0 \\ 0 & 0 & 1 
\end{array} \right) = \left. \left( \begin{array}{ccc}   
\frac{\partial^3 F_1}{\partial^2 t_1 \partial u} & \frac{\partial^3 F_1}{\partial^2 t_1 \partial x_1} & \frac{\partial^3 F_1}{\partial^2 t_1 \partial x_2}  \\
\frac{\partial^2 F_1}{\partial t_1 \partial u} & \frac{\partial^2 F_1}{\partial t_1 \partial x_1} & \frac{\partial^2 F_1}{\partial t_1 \partial x_2} \\ 
\frac{\partial F_1}{\partial u} - \frac{\partial F_2}{\partial u} & \frac{\partial F_1}{\partial x_1} - \frac{ \partial F_2}{ \partial x_1} & \frac{\partial F_1}{\partial x_2} - \frac{\partial F_2}{ \partial x_2} \\
\end{array} \right) \right|_{(A(t_i, {\bf 0}), \bx_0)} \times JB 
 $$



Now, 
$$\frac{\partial F_i}{\partial x_1} = 2X_i - 2x_1, \frac{\partial F_i}{\partial x_2} = -2Y_i + 2 x_2$$
and $\gamma_i - \bx = (X-x_1, Y-x_2) = \frac{1}{\kappa} N$ where $N = (-1)^{\beta}(Y_1', X_1')$.  Hence, 
$\frac{\partial F_i}{\partial x_1} = 2Y_1 (-1)^{\beta}$ and $\frac{\partial F_i}{\partial x_2} = -2X_1' (-1^{\beta})$. 
Substituting these derivatives into the matrix equation gives:
$$\left( \begin{array}{ccc}   1 & 0 & 0  \\ 0 & 1 & 0 \\ 0 & 0 & 1
\end{array} \right) = \left. \left( \begin{array}{ccc}  * & 2\alpha_1^2 X_1'' + 4 \alpha_2 X_1' & -2\alpha_1^2 Y_1''-4\alpha_2 Y_1'   \\ 
* & 2\alpha_1 X_1' & -2\alpha_1 Y_1'  \\ 
* & \frac{2}{\kappa}(-1)^{\beta}(Y'_1 -Y'_2) & \frac{2}{\kappa}(-1)^{\beta}(X'_2 -X'_1)
\end{array} \right) \right|_{(A(t_i, {\bf 0}), \bx_0)} \times JB. 
 $$

Recall that the type of bifurcation depends upon whether $\frac{\partial B_1}{\partial a_1} \frac{\partial B_1}{\partial a_3}$ is positive or negative.
\begin{eqnarray*}
\frac{\partial B_1}{\partial a_1} &=& \ {\rm det}\left| \begin{array}{cc}   2\alpha_1 X_1' & -2\alpha_1 Y_1'  \\ 
 \frac{2}{\kappa}(-1)^\beta(Y'_1 -Y'_2) & \frac{2}{\kappa}(-1)^\beta(X'_2 -X'_1) \end{array} \right|  \\ 
&=& 2\alpha_1 X_1' \frac{2}{\kappa}(X'_2 -X'_1) + 2\alpha_1 Y_1' \frac{2}{\kappa}(Y'_1 -Y'_2) \\
  &=& \frac{4\alpha_1}{\kappa}(-X_1'^2 + Y_1'^2   + X_1'X_2' - Y_1'Y_2')  \\ &=& \frac{4\alpha_1}{\kappa}(\langle T_1, T_1 \rangle - \langle T_1 , T_2 \rangle).
\end{eqnarray*}

and $\frac{\partial B_1}{\partial a_3} = -(2\alpha_1^2 X_1'' + 4 \alpha_2 X_1')2\alpha_1 Y_1'  + 2\alpha_1 X_1'(2\alpha_1^2 Y_1''+4\alpha_2 Y_1')=  4 \alpha_1^3 (X_1'Y_1'' - X_1'' Y_1') = 4 \alpha_1^3 \kappa.$
$$\frac{\partial B_1}{\partial a_1} \frac{\partial B_1}{\partial a_3} = 16 \alpha_1^4 (-1)^{\beta}(\langle T_1, T_1 \rangle - \langle T_1 , T_2 \rangle).$$

So if $\gamma_1$ and $\gamma_2$ are both spacelike, this gives $16 \alpha_1^4( 1 - \langle T_1 , T_2 \rangle)$ 
which is negative if $\gamma_1$ and $\gamma_2$ lie on the same branch and positive if they lie on opposite branches (see Section \ref{sectionA12A2}).  On the other hand,
if they are both timelike this gives $-16 \alpha_1^4( -1 - \langle T_1 , T_2 \rangle)$.  Parametrising the psuedo-circle of type $H^1(c, -r)$ as $H^1(\theta) = c + r(\pm \sinh(\theta), \cosh(\theta))$, the unit tangent vector is given by $T=(\pm \cosh(\theta), \sinh(\theta))$.  Now if $\gamma_1$ and $\gamma_2$ lie on the same branch
$\langle T_1 , T_2 \rangle = -\cosh(\theta_1)\cosh(\theta_2) + \sinh(\theta_1)\sinh(\theta_2) = -\cosh(\theta_1 - \theta_2)$ which is less than -1.  However if $\gamma_1$ and $\gamma_2$ lie on opposite branches
$\langle T_1 , T_2 \rangle = \cosh(\theta_1)\cosh(\theta_2) + \sinh(\theta_1)\sinh(\theta_2) = \cosh(\theta_1 + \theta_2)$ which is greater than 1.  Hence the expression $-16 \alpha_1^4( -1 - \langle T_1 , T_2 \rangle)$ is negative if $\gamma_1$ and $\gamma_2$ lie on the same branch and positive if they lie on opposite branches (the same conditions as for spacelike). It follows that the type is $A_1A_3(a)$ if both contact points lie on opposite branches and type $A_1A_3(b)$ occur on the same branch of the pseudo-circle.  


\section{The $A_4$ singularity}

Consider the following standard versal unfolding of an $A_4$ singularity given by
$$G: \mathbb{R} \times \mathbb{R}^3 \to \R $$ 
where $\mathbb{R}$ denotes the parameters $t$ and $\mathbb{R}^3$ denotes the unfolding parameters $\ba = (a_1, a_2, a_3)$ and the versal unfolding is given by the two unfoldings:

\begin{eqnarray*}
G(t, \ba) &=& t^5 +a_1 t^3 + a_2 t^2 + a_3 t.
\end{eqnarray*}

\subsection{The big bifurcation set}
\begin{figure}[tp]
\begin{center}
\includegraphics[width=5cm, height=4.1cm]{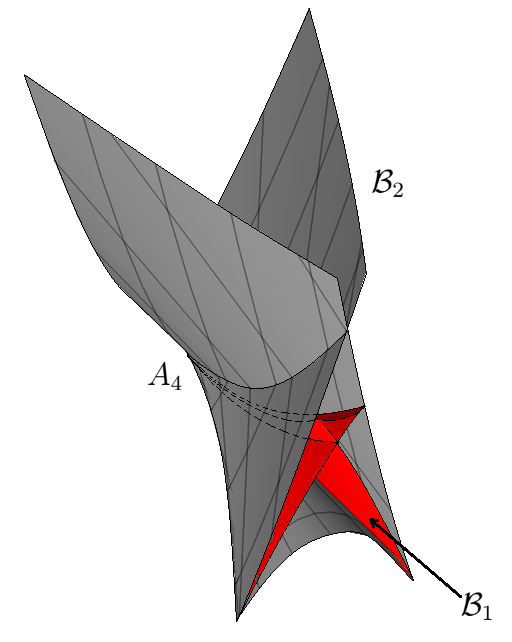}
\caption{The BBS for $A_4$.}
\label{fig:A4-BBS}
\end{center}
\end{figure}

The bifurcation set $\mathcal B_1$ of the standard $A_4$ singularity $G$ is the swallowtail surface which can be parametrised by $(a_1, -10t^3 - 3a_1t)$, and its bifurcation set $\mathcal B_2$ is another swallowtail, which sits inside the swallowtail $\mathcal B_1$ and can be parametrised by $(-3s^2 - 4st -3t^2, 2s^3 + 8s^2 t+ 8st^2 + 2t^3, -4s^3t - 7s^2t^2 - 4st^3)$. See Figure \ref{fig:A4-BBS}.
The adjacent 1-dimensional strata are found to be
\begin{eqnarray*}
  A_3  &:&  \{ (a_1, a_2, a_3) =  (-10 t^2, 20 t^3, -15 t^4) \} \\
		A_1 A_2 &:& \{ (a_1, a_2, a_3) =   (-60 t^2, -80 t^3, 960 t^4) \}   \\  
		A_2  / A_2  &:&  \{ (a_1, a_2, a_3) =  (-\frac{10}{3} t^2, 0, 5 t^4)    \} \\
	A_1^2 / A_1^2 &:& \{ (a_1, a_2, a_3) = (-4 t^2, 0, \frac{16}{5}t^4)  \}.   
\end{eqnarray*}

The limiting tangent vectors to these one-dimensional strata are all given by 
$(1, 0, 0)$, so the only bad planes is given by
$a_1 = 0$.  Examining representations from both components show that only transition type exists for $A_4$.

\end{document}